\title{Isomorphic copies of $l^\infty$ in Ces\`aro-Orlicz function spaces}
\author{Tomasz Kiwerski, Pawe\l{} Kolwicz}
\documentclass[11pt]{amsart}
\usepackage[latin1]{inputenc}
\usepackage[english]{babel}
\usepackage{amsmath,amssymb,amsthm,amsfonts}
\usepackage[all]{xy}

\makeatletter
\@namedef{subjclassname@2010}{ \textup{2010} Mathematics Subject Classification}
\makeatother

\theoremstyle{definition}

\numberwithin{equation}{section}

\numberwithin{equation}{section}

\newcommand{\norm}[1]{\left\lVert#1\right\rVert}

\begin{document}
\begin{center}
\textbf{\large{ISOMORPHIC COPIES OF $l^\infty$ IN CES\`ARO-ORLICZ\\
FUNCTION SPACES}}
\end{center}

\bigskip

\begin{center}
TOMASZ KIWERSKI
\end{center}

\begin{center}
Faculty of Mathematics, Computer Science and Econometrics\\
University of Zielona G\'{o}ra, \\
prof. Z. Szafrana 4a, 65-516 Zielona G\'{o}ra, Poland\\
e-mail address: \verb+tomasz.kiwerski@gmail.com+
\end{center}

\bigskip

\begin{center}
PAWE\L{} KOLWICZ
\end{center}

\begin{center}
Institute of Mathematics,\\
Faculty of Electrical Engineering,\\
Pozna\'{n} University of Technology,\\
Piotrowo 3A, 60-965 Pozna\'{n}, Poland\\
e-mail address: \verb+pawel.kolwicz@put.poznan.pl+
\end{center}

\bigskip
\noindent
ABSTRACT. We characterize Ces\`aro-Orlicz function spaces $Ces_\varphi$ containing isomorphic copy of $l^\infty$. We also describe the subspaces $(Ces_\varphi)_a$ of all order continuous elements of $Ces_\varphi$. Finally, we study the monotonicity structure of the spaces $Ces_\varphi$ and $(Ces_\varphi)_a$.

\bigskip

\begin{flushleft}
\end{flushleft}

\begin{flushleft}
\end{flushleft}
\bigskip

\section{INTRODUCTION}

	The structure of different types of Ces\`aro spaces has been widely investigated during the last decades from the isomorphic as well as isometric point of view. Firstly, the Ces\`aro sequence spaces $ces_p$ has been studied by many authors starting from 1970. Then the attention has been moved to the Ces\`aro function spaces $Ces_p$. It was interesting among others that some properties are fulfilled in the sequence case and do not in the function case. Moreover, it happens that the cases $Ces_p[0,1]$ and $Ces_p[0,\infty)$ are essentially different (see the description of the K\"{o}the dual of Ces\`aro spaces in \cite{AM09} and \cite{LM15a}). Probably the most important papers concerning the structure of the spaces $Ces_p$ are two papers by Astashkin and Maligranda (see \cite{AM09}, \cite{AM14}). Recall that the space $Ces_p[0,1]$, for $1\le p \le \infty$, consists of those Lebesgue measurable real functions $f$ on $[0,1]$ for which the Ces\`aro-means $C|f|(x)=\frac{1}{x} \int_0^x|f(t)|dt$ belong to $L_p[0,1]$ (similarly for the spaces $ces_p$ and $Ces_p[0,\infty)$). It is natural to investigated the space
\begin{equation*}
CX=CX(I)=\{f \in L^0(I) : C|f| \in X\}
\end{equation*}
for any Banach ideal space $X$ on $I$, where  $I=[0,1]$  or $I=[0,\infty)$. These spaces have been defined in \cite{Ru80} for $X$ being a Banach ideal space on $[0,\infty)$ and, for example, in \cite{DS07}, \cite{LM15a} for $X$ being a symmetric space on $I$.  If we take instead of $X$ the Orlicz sequence space $l^\varphi$ (the generalization of $l^p$) then 
the space $CX$ becomes the Ces\`aro-Orlicz sequence space denoted by $ces_\varphi$ which the structure is of 
course more rich than for the Ces\`aro sequence space $ces_p$. The spaces $ces_\varphi$ has been studied by many authors (e.g. \cite{CHPSSz05}, \cite{CJP97}, \cite{CMP00}, \cite{KK10} and \cite{Ku09}). If we put $X=L^\varphi(I)$ the Orlicz function space then $CX$ becomes the Ces\`aro-Orlicz function space $Ces_\varphi(I)$. As far as we know the structure of these spaces has 
not been investigated until now. We want to characterize these spaces $Ces_\varphi$ which contain an order isomorphic copy of $l^\infty$ (equivalently, are not order continuous). Clearly, such results are very useful in further studying isomorphic structure of these spaces. We characterize the subspaces $(Ces_\varphi)_a$ of all order continuous elements in $Ces_\varphi$. Finally, we prove criteria for strict monotonicity of the spaces $(Ces_\varphi)_a$. We also give a characterization of uniform monotonicity in Ces\`aro-Orlicz function spaces. We consider the largest possible class of Orlicz functions giving the maximal generality of spaces under consideration.

\section{PRELIMINARIES}

The symbols $\mathbb{R}$, $\mathbb{R_+}$ and $\mathbb{N}$ denote the sets of reals, nonnegative reals and natural numbers, respectively.
Let $L^0=L^0 (I)$ be the space of all classes of real-valued Lebesgue measurable functions defined on $I$, where 
$I=[0,1]$ or $I=[0,\infty)$. A Banach space $E=(E,\|\cdot\|)$ is
said to be a Banach ideal space on $I$ if $E$ is a linear subspace of $L^0(I)$ and satisfies two conditions:
\begin{enumerate} 
  \item if $g\in E$, $f\in L^0$ and $|f|\leq |g|$ a.e. on $I$ then $f \in E$ and $\|f\|\leq \|g\|$,
  \item there is an element $f \in E$ that is positive on whole $I$ .
\end{enumerate}
Sometimes we write $\norm{\cdot}_{E}$ to be sure in which space the norm has been taken.

For two Banach ideal spaces $E$ and $F$ on $I$ the symbol $E \hookrightarrow F$ means that the embedding $E\subset F$ is continuous, i.e., there exists constant a $C>0$ such that $\norm{x}_{F} \le C\norm{x}_{E}$ for all $x\in E$. Moreover, $E=F$ means that the spaces are the same as the sets and the norms are equivalent.

Recall that for $f\in E$ the distribution function $d_f$ is defined by
\begin{equation*}
d_f(\lambda)=m(\{t\in I : |f(t)|>\lambda\})
\end{equation*}
for all $\lambda>0$, where $m$ is the Lebesgue measure. The non-increasing rearrangement of $f$ is denoted by $f^*$ and is defined as
\begin{equation*}
f^*(t)=\inf\{\lambda >0 : d_f(\lambda)<t\}
\end{equation*}
for $t \ge 0$. We say that two functions $f,g \in L^0(I)$ are equimeasurable if they have the same distribution functions $d_f \equiv d_g$. Then we write $f\sim g$. By a symmetric function space (symmetric Banach function space or rearrangement invariant Banach function space) on $I$ we mean a Banach ideal space $(E,\norm{\cdot}_{E})$ with an additional property that for any two functions $f\in E$,$g\in L^0(I)$, with $f\sim g$ we have $g\in E$ and $\norm{f}_{E} = \norm{g}_{E}$. In particular, $\norm{f}_{E} = \norm{f^*}_{E}$.

A point $f\in E$ is said to have an order continuous norm (or to be an order continuous element) if for each sequence $(f_n)\subset E$ satisfying $0\le f_n\le |f|$ and $f_n\rightarrow 0$ a.e. on I, one has $\norm{f_n}\rightarrow 0$. By $E_a$ we denote the subspace of all order continuous elements of $E$. It is worth to notice that in case of Banach ideal spaces on $I$, $x\in E_a$ if and only if $\norm{x\chi_{A_n}}\rightarrow 0$ for any decreasing sequence of Lebesgue measurable sets $A_n\subset I$ with empty intersection. A Banach ideal space is called order continuous ($E\in \text{(OC)}$ shortly) if every element of $E$ is order continuous, i.e., $E=E_a$.

We say a Banach ideal space $X$ is strictly monotone $(X\in \text{(SM)})$, if $\norm{x} < \norm{y}$ for all $x,y\in X$ such that$0\le  x\le y$ and $x\neq y$. A Banach ideal space $X$ is uniformly monotone $(X\in \text{(UM)})$ if for each $\epsilon \in (0,1)$ there exists $\delta(\epsilon)\in (0,1)$ such that for all $0\le y\le x$, $\norm{x}=1$ and $\norm{y}\ge \epsilon$, we have $\norm{x-y}\le 1-\delta(\epsilon)$. $X$ is called lower locally uniformly monotone $(X\in \text{(LLUM)})$ if for any $x\in X$, $\norm{x}=1$ and $\epsilon \in (0,1)$ there exists $\delta(x,\epsilon) \in (0,1)$) such that $\norm{x-y}\le 1-\delta(x,\epsilon)$ whenever $0\le y\le x$ and $\norm{y}\ge \epsilon$. Note that monotone properties of Banach ideal spaces are useful in dominated best approximation problems, see \cite{CKP14} for further references.

The next well known theorem shows the relation between the order continuity of E and the existing of isomorphic copy of $l^\infty$.
\bigskip
\\\textbf{Theorem A}. (G. Ya. Lozanovski{\u \i}, see \cite{Lo69}) A Banach ideal space $E$ is order continuous if and only if $E$ contains no isomorphic copy of $l^\infty$.
\bigskip
\\A function $\varphi : [0,\infty)\rightarrow[0,\infty]$ is called an Orlicz function if:
\begin{enumerate}
\item $\varphi$ is convex,
\item $\varphi$ is non-decreasing,
\item $\varphi(0)=0$,
\item $\varphi$ is neither identically equal to zero nor infinity on $(0,\infty)$,
\item $\varphi$ is left continuous on $(0,\infty)$, i.e., $\lim_{u\rightarrow b^{-}_{\varphi}} \varphi(u)=\varphi(b_\varphi)$
if $b_\varphi <\infty$, where
\begin{equation*}
b_\varphi = \sup\{u>0 : \varphi(u)<\infty\}.
\end{equation*}
\end{enumerate}
For more information about Orlicz functions see \cite{Ch96} and \cite{KR61}.

If we denote
\begin{equation*}
a_\varphi=\sup\{u\ge 0 : \varphi(u)=0\},
\end{equation*}
then $0\le a_\varphi\le b_\varphi\le \infty$, $a_\varphi<\infty$, $b_\varphi>0$, since an Orlicz function is neither identically equal to zero nor infinity on $(0,\infty)$. The function $\varphi$ is continuous and nondecreasing on $[0,b_\varphi)$ and is strictly increasing on $[a_\varphi,b_\varphi)$. We use notations $\varphi>0$, $\varphi<\infty$ when $a_\varphi=0$, $b_\varphi<\infty$, respectively.

We say an Orlicz function $\varphi$ satisfies the condition $\Delta_2$ for large arguments $(\varphi\in \Delta_2(\infty) \text{for short})$ if there exists $K>0$ and $u_0>0$ such that $\varphi(u_0 )<\infty$ and
\begin{equation*}
\varphi(2u)\le K\varphi(u)
\end{equation*}
for all $u\in [u_0,\infty)$. Similarly, we can define the condition $\Delta_2$ for small, with $\varphi(u_0)>0$ $(\varphi\in \Delta_2 (0))$ or for all arguments $(\varphi \in \Delta_2 (\mathbb{R_+}))$. These conditions play a crucial role in the theory of Orlicz spaces, see \cite{Ch96}, \cite{KR61}, \cite{Ma89} and \cite{Mu83}.

We will write $\varphi \in \Delta_2$ in two cases: $\varphi\in \Delta_2(\infty)$ if $I=[0,1]$ and $\varphi \in \Delta_2 (\mathbb{R_+})$ if $I=[0,\infty)$.

The Orlicz function space  $L^\varphi=L^\varphi(I)$ generated by an Orlicz function $\varphi$ is defined by
\begin{equation*}
L^\varphi = \{f\in L^0(I) : I_\varphi(f/\lambda)<\infty \ \text{for some}\  \lambda=\lambda(f)>0\},
\end{equation*}
where $I_\varphi(f)=\int_I \varphi(|f(t)|)dt$ is a convex modular (for the theory of Orlicz spaces and modular spaces see \cite{Ma89} and \cite{Mu83}). The space $L^\varphi$ is a Banach ideal space with the Luxemburg-Nakano norm
\begin{equation*}
\norm{f}_{\varphi}=\inf\{\lambda : I_\varphi(f/\lambda)\le 1\}.
\end{equation*}
It is well known that $\norm{f}_{\varphi}\le 1$ if and only if $I_\varphi(f)\le 1$. Moreover, the set
\begin{equation*}
KL^\varphi=KL^\varphi(I)=\{f\in L^0(I) : I_\varphi(f)<\infty\},
\end{equation*}
will be called the Orlicz class.

The continuous Ces\`aro operator is defined for $0<x\in I$ as
\begin{equation*}
Cf(x)=\frac{1}{x}\int_0^x f(t)dt.
\end{equation*}
For a Banach ideal space $X$ on $I$ we define an abstract Ces\`aro space $CX=CX(I)$ as
\begin{equation*}
CX=\{f\in L^0(I) : C|f|\in X\}
\end{equation*}
with a norm $\norm{f}_{CX}=\norm{C|f|}_{X}$.

The Ces\`aro-Orlicz function space $Ces_\varphi=Ces_\varphi(I)$ is defined as $Ces_\varphi(I)=CL^\varphi(I)$. Consequently, the norm in the space $Ces_\varphi$ is given by the formula
\begin{equation*}
\norm{f}_{Ces(\varphi)}=\inf\{\lambda : \rho_\varphi(f/\lambda)\le 1\}
\end{equation*}
where $\rho_\varphi(f)=I_\varphi(C|f|)$ is a convex modular.

The dilation operator $D_s$, $s>0$, defined on $L^0 (I)$ by
\begin{equation*}
D_s x(t)=x(t/s) \chi_I (t/s)=x(t/s) \chi_{[0,\min\{1,s\} )}  (t),
\end{equation*}
for $t\in I$, is bounded in any symmetric space $E$ on $I$ and $\norm{D_s}_{E\rightarrow E}\le \text{max}\{1,s\}$ (see \cite[p. 148]{BS88}). Moreover, the lower and upper Boyd indices of $E$ are defined by
\begin{equation*}
p(E)=\lim \limits_{s\rightarrow 0^{+}} \frac{\text{ln} \norm{D_s}_{E\rightarrow E}}{\text{ln}s},
\end{equation*}
\begin{equation*}
q(E)=\lim \limits_{s\rightarrow \infty} \frac{\text{ln} \norm{D_s}_{E\rightarrow E}}{\text{ln}s}.
\end{equation*}
In particular, they satisfy the inequalities
\begin{equation*}
1\le p(E)\le q(E) \le \infty.
\end{equation*}
In the case when $E$ is the Orlicz space $L^\varphi$, these indices correspond to the so-called Orlicz-Matuszewska indices of Orlicz functions generating the Orlicz spaces, i.e., $\alpha_\varphi=p(L^\varphi)$ and $\beta_\varphi=q(L^\varphi)$, where $\alpha_\varphi$ and $\beta_\varphi$ are the lower and upper Orlicz-Matuszewska indices of the Orlicz space $L^\varphi$ (see \cite[Proposition 2.b.5 and Remark 2 on page 140]{LT79}). For more details see \cite{Bo69}, \cite{Bo71}, \cite{Ma85}, \cite{Ma89} and \cite{Mu83}.

In this paper we accept the convention that $\sum_{n=m}^{k} x_n=0$  if $k<m$.
\bigskip
\\Let us mention the important result about boundedness of the operator $C$.
\bigskip
\\\textbf{Theorem B.} \cite[\text{p.} 127]{KMP07} For any symmetric space $E$ on $I$ the operator $C : E\rightarrow E$ is bounded if and only if the lower Boyd index satisfies $p(E)>1$.
\bigskip
\\The immediate consequence of Theorem $B$ and the above discussion is a next corollary.
\bigskip
\\\textbf{Corollary 1}. The embedding $L^\varphi \hookrightarrow Ces_\varphi$ holds if and only if $\alpha_\varphi>1$.
\bigskip
\\\textbf{Remark.} Let $X$,$Y$ be Banach ideal spaces on $I$. If $X\hookrightarrow Y$ then $CX\hookrightarrow CY$. Indeed, suppose  $X\hookrightarrow Y$. Then for all $x\in CX$
\begin{equation*}
\norm{x}_{CY}=\norm{C|x|}_{Y}\le A\norm{C|x|}_{X}=A\norm{x}_{CX},
\end{equation*}
for some constant $A>0$. This means that $CX\hookrightarrow CY$.

Let $L^\varphi$ and $L^\psi$ be the Orlicz spaces with $\varphi,\psi<\infty$. The criteria for the embeddings $L^\varphi\hookrightarrow L^\psi$ can be found in \cite[Theorem 3.4 p. 18]{Ma89}. Consequently, we conclude that:
\begin{enumerate}
\item if there exists $k>0$ with $\psi(u)\le \varphi(ku)$ for all $u\in [0,\infty)$ then
$Ces_\varphi [0,\infty)\hookrightarrow Ces_\psi [0,\infty)$,
\item if there exists $k,u_0>0$ such that $\psi(u)\le \varphi(ku)$ for all $u\in [u_0,\infty)$ then 
$Ces_\varphi [0,1]\hookrightarrow Ces_\psi [0,1]$,
\item if there exists $k,u_0>0$ such that $\psi(u)\le \varphi(ku)$ for all $u\in [0,u_0 ]$ then 
$ces_\varphi\hookrightarrow ces_\psi$.
\end{enumerate}
Note that case (iii) is exactly Proposition 1 in \cite{Ku09} but now it is an immediate consequence of our above remark.

\section{RESULT}

\noindent
\textbf{Fact.} Let $X$ be a Banach ideal space on $I$. By the definition, the order continuity of $X$ gives the same for $CX$ (for proof see [19, Lemma 1]). The case of strict monotonicity and uniform monotonicity is similar. The converse is not true in general (see [11, Proposition 2.1] and [21, Example 1]).
\begin{proof}
We proof only the implication for uniform monotonicity. We apply Theorem 6 (ii) in [13]. Suppose that $X\in \text{(UM)}$. Take $\epsilon>0$, $x,y\in CX$, $x,y\ge 0$, $\norm{x}_{CX}=1$ and $\norm{y}_{CX}\ge \epsilon$. We have
\begin{equation*}
\norm{x+y}_{CX}=\norm{C|x+y|}_{X}=\norm{C|x|+C|y|}_{X}.
\end{equation*}
Since $X$ is uniformly monotone we have that there is a  $\sigma(\epsilon)>0$ such that
\begin{equation*}
\norm{C|x|+C|y|}_{X}\ge 1+\sigma(\epsilon).
\end{equation*}
But this means that $CX\in \text{(UM)}$. 
\end{proof}
\bigskip
\noindent
\textbf{Proposition 2.} Suppose $X$ is a symmetric space on $I=[0,1]$ and $C : X\rightarrow X$. Then $X\in \text{(OC)}$ if and only if $CX\in \text{(OC)}$.
\begin{proof}
\emph{Necessity}. It follows from Fact above.

\noindent
\emph{Sufficiency}. We thank Professor Karol Le\'{s}nik for giving the proof.

\noindent
Suppose that $X\notin \text{(OC)}$, i.e. there exists $0< f\in X\backslash X_a$. Because $X$ is symmetric we can assume that $f=f^*$ (see [7, Lemma 2.6]). Then $Cf\in X$ and $Cf\ge f$. Applying the symmetry of $X$ and passing to subsequence, if necessary, we can assume that there is a $\delta>0$ such that
$$\norm{f\chi_{[0,1/n)}}_{X}\ge \delta,$$
for all $n\in \mathbb{N}$. Since
$$C(f\chi_{[0,1/n)})(t)\ge f\chi_{[0,1/n)}(t),$$
for all $t\in I$ and $n\in \mathbb{N}$, we have
$$\norm{f\chi_{[0,1/n)}}_{CX}=\norm{C(f\chi_{[0,1/n)})}_{X}\ge \norm{f\chi_{[0,1/n)}}_{X}\ge \delta >0,$$
which means that $X\notin \text{(OC)}$.                  
\end{proof}

\noindent
\textbf{Proposition 3.} Let $\varphi$ be an Orlicz function. The following conditions are equivalent:
\begin{enumerate}
\item the space $Ces_\varphi [0,\infty)\neq{0}$,
\item there exists $\lambda_0>0$ and $x_0\in [0,\infty)$ such that $\int_{x_0}^\infty \varphi(\lambda_0/t)dt<\infty$,
\item for each $\lambda_0>0$ there exist $y_0\in [0,\infty)$ with $\int_{y_0}^\infty \varphi(\lambda_0/t)dt<\infty$.
\end{enumerate}
\begin{proof}
The equivalence of conditions (i) and (ii) follows from Theorem 1 $(a)$ in [19]. Clearly, $\text{(iii)}\Rightarrow \text{(ii)}$.

\noindent
$\text{(ii)}\Rightarrow \text{(iii)}$. We have to consider two cases. If $\lambda\le \lambda_0$, there is nothing to proof. Now suppose $\lambda> \lambda_0$. Then
\begin{equation*}
\int_{y_0}^\infty \varphi(\lambda/t)dt=\int_{y_0}^\infty \varphi\left(\frac{\lambda_0}{\lambda_0t/\lambda}\right)dt=\frac{\lambda}{\lambda_0}
\int_{\lambda_0 y_0/\lambda}^{\infty}\varphi(\lambda_0/u)du
\end{equation*}
therefore it is enough to take $y_0=(\lambda/\lambda_0)x_0$.
\end{proof}

\noindent
\textbf{Remark 3.1.} Let $\varphi$ be an Orlicz function satisfying condition $(S)$:
$$\liminf \limits_{t\rightarrow 0} t\varphi^{'}(t)/\varphi(t)>1.$$
Then $Ces_\varphi [0,\infty)\neq{0}$.
\begin{proof}
Firstly, condition $(S)$ implies (ii) in Proposition 3. Indeed, we can use the same argument as in the proof of implication $\text{(a)}\Rightarrow\text{(c)}$ in [8, Theorem 2.2]. But condition (ii) is equivalent to $Ces_\varphi [0,\infty)\neq{0}$.	
\end{proof}

\noindent
\textbf{Remark 4.} Let $\varphi$ be an Orlicz function. Then $Ces_\varphi [0,1]\neq{0}$.
\begin{proof}
It follows from Theorem 1 (b) in [19] that $Ces_\varphi [0,1]\neq{0}$ if and only if there exists $0<a<1$ such that  $\chi_{[a,1]} \in L^\varphi [0,1]$. But $L^\varphi [0,1]$ is a symmetric space and
$$L^\infty [0,1]=L^1 [0,1]\cap L^\infty [0,1]\hookrightarrow L^\varphi [0,1]\hookrightarrow L^1 [0,1]+L^\infty [0,1]=L^1 [0,1]$$
(see [3]). Therefore $\chi_{[a,1]} \in L^\varphi [0,1]$ for all $a\in (0,1)$.
\end{proof}

\noindent
The following theorem gives a similar information as in the sequence space $ces_\varphi$ one can get from Theorem 2.3 in [8]. Before we formulate the theorem define a set
$$C_\varphi=C_\varphi (I)=\{x\in Ces_\varphi (I) : \rho_\varphi (kx)<\infty \ \text{for all}\  k>0\}.$$
Note that $C_\varphi= \{0\}$ whenever $b_\varphi <\infty $.

\bigskip
\noindent
\textbf{Theorem 5.} Suppose $\varphi<\infty$. Then the following conditions are true:
\begin{enumerate}
\item $C_\varphi$ is the subspace of all order continuous elements of $Ces_\varphi$, 
\item $C_\varphi$ is a closed separable subspace of $Ces_\varphi$.
\end{enumerate}
\begin{proof}
(i). We will show that $C_\varphi \hookrightarrow (Ces_\varphi)_a$. Take any $x\in C_\varphi$ and a sequence $(A_n )$ with $A_n\searrow \emptyset$. It is enough to show, that $\norm{x\chi_{A_n}}_{Ces(\varphi)}\rightarrow 0$ or equivalently $\rho_\varphi (kx \chi_{A_n} )\rightarrow 0$ for any $k>0$. Set $k>0$. Clearly,
$$\lim \limits_{n\rightarrow \infty} m((0,t]\cap A_n )=0,$$
for each $t>0$. Then for almost each $t>0$
$$0\le C(|x\chi_{A_n} ) |)(t)=\frac{1}{t} \int_0^t|x(s) \chi_{A_n}(s) |ds \rightarrow 0.$$
Moreover,
$$\rho_\varphi (kx\chi_{A_n})\le \rho_\varphi (kx)<\infty,$$
for each $k$. Thus, by the  Lebesgue dominated convergence theorem, we have $\rho_\varphi (kx\chi_{A_n})\rightarrow 0$, so $x\in (Ces_\varphi (I))_a$.

\noindent
Next we prove the reverse inclusion. Define a set
$$B_\varphi=B_\varphi(I)=\text{cl}\{f\in Ces_\varphi(I) : f \ \text{is simple function and} \ m(\text{supp}(f))<\infty\}.$$
We claim that :
$$ C_\varphi=B_\varphi. \hspace{1.5cm} (A)$$

Consider the inclusion $B_\varphi \subset C_\varphi$. Of course, if $B_\varphi=\emptyset$ then $B_\varphi \subset C_\varphi$ so we can assume that $B_\varphi \neq \emptyset$. We divide the proof in three cases.
\begin{enumerate}
\item[(1)] Suppose $x=\chi_[a,b]$ , where $[a,b]\subset I$, $0\le a<b<\infty$. Since $\varphi < \infty$ and $Ces_\varphi \neq {0}$, by Proposition 3 we conclude that $\rho_\varphi (kx)<\infty$ for all $k>0$.
\item[(2)] Let $x=\sum_{k=1}^{n} c_k \chi_{A_k}$, $A_k=[a_k,b_k ]\subset I$ for $n,k\in \mathbb{N}$ and $0\le a_k<b_k<\infty$. Using the convexity of modular $\rho_\varphi$ and argument  from case (1), we obtain $\rho_\varphi (\lambda x)<\infty$ for all $\lambda >0$.
\item[(3)] Assume that there exists a sequence $(x_n )$ such that  $ x_n\rightarrow x$ in $Ces_\varphi$ and $x_n$ are of the form (2). Let $k>0$. For sufficiency large $n\in \mathbb{N}$ we have
\begin{equation*}
\rho_\varphi(kx)=\rho_\varphi(k(x-x_n+x_n))=\rho_\varphi\left(\frac{1}{2}2k(x-x_n )+\frac{1}{2}2kx_n\right)
\end{equation*}
\begin{equation*}
\le \frac{1}{2}\rho_\varphi(2k(x-x_n ))+\frac{1}{2}\rho_\varphi(2kx_n)\le 1+\rho_\varphi(2kx_n)<\infty.
\end{equation*}
Therefore $B_\varphi \subset C_\varphi$.
\end{enumerate}

Now, we will prove the reverse inclusion $C_\varphi \subset B_\varphi$. Let $x\in C_\varphi$, $k>0$ and $\epsilon>0$. Since $\int_I \varphi(C(k|x|))dt<\infty$ there exists $\beta=\beta(k)\in \mathbb{R_+}$ with
$$\int_{\beta}^\infty \varphi(C(2k|x(t)|))dt<\epsilon/2.$$
Put $z_n=\varphi(C(2k|x|\chi_{[0,1/n)}))$ and $z=\varphi(C(2k|x|))$. Then $z\in L^1$, $0\le z_n \le z$ and $z_n \rightarrow 0$ a.e. on $I$. Note that $L^1\in \text{(OC)}$ whence $\norm{z_n}_{L^1}\rightarrow 0$. Consequently, there exists $\alpha=\alpha(k)\in \mathbb{R_+}$ such that
$$\int_{\beta}^\infty \varphi(C(2k|x|\chi_{[0,\alpha)}))dt<\epsilon/2.$$
Define an element $x_n=x\chi_{[1/n,n]}$ for $n\in \mathbb{N}$. For sufficiently large $n\in \mathbb{N}$ satisfying $1/n<\alpha$ and $n>\beta$ we have
$$\rho_\varphi(k(x_n-x))=\rho_\varphi(kx\chi_{[0,1/n)\cup (n,\infty)})=\int_0^\infty \varphi(C(kx(t)\chi_{[0,1/n)\cup (n,\infty)}(t)))dt$$
$$=\int_{0}^{\infty} \varphi(C(kx(t)\chi_{[0,1/n)}(t))+C(kx(t)\chi_{(n,\infty)}(t)))dt$$
$$\le \frac{1}{2}\int_{0}^{\infty} \varphi(C(kx(t)\chi_{[0,1/n)}(t)))dt+\frac{1}{2}\int_{0}^{\infty} \varphi(C(kx(t)\chi_{(n,\infty)}(t)))dt$$
$$\le \frac{1}{2}\int_{0}^{\infty} \varphi(C(kx(t)\chi_{[0,\alpha)}(t)))dt+\frac{1}{2}\int_{0}^{\infty} \varphi(C(kx(t)\chi_{(\beta,\infty)}(t)))dt<\epsilon$$
i.e. $\norm{x_n-x}_{Ces(\varphi)}\rightarrow 0$. Now, it is enough to prove that for each interval $[a,b]$, $0\le a<b<\infty$, there is a sequence of functions $y_n$ of the form (2) with
$$\norm{x\chi_{[a,b]}-y_n}_{Ces(\varphi)}\rightarrow 0,$$
as $n\rightarrow 0$. Without loss of generality we may assume that $x\chi_{[a,b]}\ge 0$. Since $x\in L^0$ we can find a sequence $(y_n )$ of simple functions such that $y_n\nearrow x \chi_{[a,b]}$. Thus
$$0\le x\chi_{[a,b]}-y_n\le x\chi_{[a,b]},$$
and $x\chi_{[a,b]}-y_n\rightarrow 0$ a.e. on $I$. But $x\in C_\varphi \subset (Ces_\varphi(I))_a$, whence
$$\norm{x\chi_{[a,b]}-y_n}_{Ces(\varphi)}\rightarrow 0.$$
This proves the claim (A).

Take $x\in (Ces_\varphi(I))_a $. Define an element $x_n=x\chi_{[1/n,n]}$. Then $\norm{x_n-x}_{Ces(\varphi)}\rightarrow 0$. Moreover, $x_n\in B_\varphi$ for all $n\in \mathbb{N}$ by the above proof. Since $B_\varphi$ is closed, so $x\in B_\varphi=C_\varphi$.

\noindent
(ii) Clearly, since $\rho_\varphi$ is a convex modular, $C_\varphi$ is a linear subspace of $Ces_\varphi$. $C_\varphi$ is closed, since $B_\varphi=C_\varphi$ and $B_\varphi$ is closed by the definition (note that $X_a$ is closed for each Banach ideal space $X$, see [3, proof of Theorem 3.8, p. 16]). Separability of the space $C_\varphi$ follows from (i) and the fact, that Lebesgue measure is separable, see [3, Theorem 5.5, p. 27]. 
\end{proof}

\noindent
\textbf{Remark 6.} Note that from Theorem 3.11, p. 18 in [3] we have the following inclusions
$$(Ces_\varphi)_a\subset (Ces_\varphi)_b\subset Ces_\varphi.$$
But $B_\varphi =(C_\varphi)_b$ (compare with the Definition 3.9, p. 17 in [3]) and $B_\varphi\subset C_\varphi\subset (Ces_\varphi)_a$ from some parts of the proof of Theorem 5. Therefore, by different argumentation we also get $C_\varphi=(Ces_\varphi)_a$.
\bigskip
\\
\noindent
\textbf{Theorem 7.} Let $\varphi$ be an Orlicz function and suppose $C : L^\varphi\rightarrow L^\varphi$. Then the following conditions are equivalent:
\begin{enumerate}
\item $\varphi \in \Delta_2$,
\item $Ces_\varphi \in \text{(OC)}$,
\item $Ces_\varphi$ contains no isomorphic copy of $l^\infty$.
\end{enumerate}
\begin{proof}
Equivalence $\text{(ii)}\Leftrightarrow \text{(iii)}$ follows from Theorem A. Now we will proof the equivalence of conditions (i) and (ii).

\noindent
$\text{(i)}\Rightarrow\text{(ii)}$. If $\varphi \in \Delta_2$, then $L^\varphi$ is order continuous, see e.g. [26, p. 21-22]. Therefore also $Ces_\varphi$ is order continuous by the Fact.

\noindent
$\text{(ii)}\Rightarrow\text{(i)}$. We have to consider two cases.

I. Suppose $I=[0,\infty)$.

\noindent
(1). Let $b_\varphi<\infty$ and $\varphi(b_\varphi)=\infty$. Let $(u_n)\subset \mathbb{R_+}$ be the sequence with $u_n\nearrow b_\varphi^{-}$. Of course, $\varphi(u_n)\rightarrow \infty$. Passing to a subsequence if necessary, we can assume that $\varphi(u_n)\ge 1$ for all $n\in \mathbb{N}$. Let $a_n=1/2^n\varphi(u_n)$ for $n\in \mathbb{N}$ and denote $a=\sum_{n=1}^\infty a_n$. Define a sequence of pairwise disjoint open intervals $(A_n )_{n\in \mathbb{N}}\subset [0,1]$ as follows
$$A_n=\left(a-\sum_{k=1}^n \frac{1}{2^k\varphi(u_k)}, a-\sum_{k=1}^{n-1} \frac{1}{2^k\varphi(u_k)}\right),$$
for all $n\in \mathbb{N}$. Then $m(A_n )=a_n$  for all $n\in \mathbb{N}$. Let
$$x=\sum_{n=1}^{\infty} u_n\chi_{A_n}.$$
We claim that $x\in Ces_\varphi$. In fact,
$$I_\varphi(x)=\int_0^\infty \varphi(x) dm=\int_0^\infty \sum_{n=1}^\infty \varphi(u_n)\chi_{A_n}dm$$
$$=\sum_{n=1}^\infty \int_{A_n} \varphi(u_n)dm=\sum_{n=1}^\infty \varphi(u_n)m(A_n)=\sum_{n=1}^\infty \frac{1}{2^n}=1.$$
The assumption $C : L^\varphi\rightarrow L^\varphi$ implies that $L^\varphi\hookrightarrow Ces_\varphi$, so $x\in Ces_\varphi$. Note that there exists $N\in \mathbb{N}$ such that for all $k\ge N$, $2u_k>b_\varphi$. Consequently, since $x$ is non-increasing
$$\rho_\varphi(2x)=I_\varphi(C(2x))=I_\varphi(2Cx)\ge I_\varphi(2x)=\sum_{n=1}^\infty \varphi(2u_n)m(A_n)=\infty,$$
and, by Theorem 5 (i), we conclude that $Ces_\varphi \notin \text{(OC)}$.

\noindent
(2). Assume that $b_\varphi<\infty$ and $\varphi(b_\varphi)<\infty$. Define an element $x=b_\varphi\chi_{[0,1]}$. Then
$$Cx(t)=b_\varphi\chi_{[0,1]}(t)+\frac{b_\varphi}{t}\chi_{(1,\infty)}(t),$$
and
$$\rho_\varphi(x)=I_\varphi(Cx)=\varphi(b_\varphi)+\int_{1}^{\infty}\varphi\left(\frac{b_\varphi}{t}\right)dt.$$
Since $I_\varphi(x)<\infty$ so $x\in L^\varphi$ and $x\in Ces_\varphi$ because $C : L^\varphi \rightarrow L^\varphi$. Moreover,
$$\rho_\varphi(2x)=I_\varphi(C(2x))\ge I_\varphi(2x)=\varphi(2b_\varphi)=\infty,$$
thus, $Ces_\varphi [0,\infty)\notin \text{(OC)}$, by Theorem 5 (i).

\noindent
(3) Let $a_\varphi>0$. Put $x=a_\varphi\chi_{[0,\infty)}$. Then $x\in Ces_\varphi [0,\infty)$ and
$$\rho_\varphi(2x)=I_\varphi(C(2x))=I_\varphi(2x)=\varphi(2a_\varphi)m([0,\infty))=\infty.$$
Therefore, $Ces_\varphi [0,\infty)\notin \text{(OC)}$ by Theorem 5 again.

\noindent
(4). Now we assume that $\varphi>0$, $\varphi<\infty$ and $\varphi \notin \Delta_2 (\mathbb{R_+})$. This means, that $\varphi \notin \Delta_2(0)$ or $\varphi \notin \Delta_2 (\infty)$. If $\varphi \notin \Delta_2 (\infty)$ then
$$\varphi(2u_n )\ge 2^n \varphi(u_n ),$$
for some sequence $(u_n )\nearrow\infty$. Therefore in this case we can use the arguments from I (1) above. Now assume that $\varphi \notin \Delta_2(0)$. That means, that there exist a decreasing sequence $(u_n )_{n\in\mathbb{N}}\subset\mathbb{R_+}$, $u_n\searrow0$ and $\varphi(2u_n )\ge 2^n \varphi(u_n )$ for all $n\in\mathbb{N}$. Define an element
$$x=\sum_{n=1}^\infty u_n \chi_{B_n},$$
where $B_n=\left(\sum_{k=1}^{n-1} \frac{1}{2^k \varphi(u_k )}, \sum_{k=1}^{n} \frac{1}{2^k \varphi(u_k )}\right)$. Then
$$I_\varphi (x)=\int_{0}^\infty \varphi(x)dm=\int_{0}^\infty \sum_{n=1}^\infty \varphi(u_n)\chi_{B_n} dm$$
$$=\sum_{n=1}^\infty \int_{B_n} \varphi(u_n)dm = \sum_{n=1}^\infty \varphi(u_n)m(B_n) = \sum_{n=1}^\infty \frac{1}{2^n} = 1,$$
so $x\in Ces_\varphi$ by the assumption $C : L^\varphi \rightarrow L^\varphi$. Furthermore, since $x$ is non-increasing,
$$\rho_\varphi(2x)=I_\varphi(C(2x))\ge I_\varphi(2x)$$
$$=\sum_{n=1}^\infty \varphi(2u_n)m(B_n) \ge \sum_{n=1}^\infty 2^n \varphi(u_n)m(B_n) = \sum_{n=1}^\infty 1 = \infty,$$
which means that the space $Ces_\varphi [0,\infty)$ is not order continuous by Theorem 5 above.

II. Suppose $I=[0,1]$.

\noindent
(1). Let $b_\varphi<\infty$ and $\varphi(b_\varphi)=\infty$. In this case we can use construction from case I (1) to show that $Ces_\varphi \in \text{(OC)}$.

\noindent
(2). Assume that $b_\varphi<\infty$ and $\varphi(b_\varphi)<\infty$. Now we follow as in the proof of case I (2) but in this case we don't need the assumption that $C : L^\varphi\rightarrow L^\varphi$. In fact, put $x=b_\varphi \chi_{[0,1]}$ . Then
$$\rho_\varphi(x)=I_\varphi(Cx)=I_\varphi(x)=\int_{0}^{1} \varphi(b_\varphi)dm=\varphi(b_\varphi)<\infty,$$
so $x\in Ces_\varphi$. Furthermore,
$$\rho_\varphi(2x)=I_\varphi(C(2x))=I_\varphi(2C(x))=I_\varphi(2x)=\varphi(2b_\varphi)=\infty,$$
which means that $Ces_\varphi\notin\text{(OC)}$.

\noindent
(3) Now suppose $\varphi<\infty$ and $\varphi \notin \Delta_2 (\infty)$. It is well known that $L^\varphi [0,1]\in \text{(OC)}$ if and only if $\varphi \in \Delta_2 (\infty)$ (see, e.g. [26, p. 21-22]). Therefore, from Proposition 2, $Ces_\varphi[0,1]\notin\text{(OC)}$. Additionally, the same direct proof of case I (4) above works also for $I=[0,1]$.   	

\end{proof}

\noindent
An immediate consequence of Fact, Theorem 5 and Theorem 7 is the next corollary.
\bigskip
\\
\noindent
\textbf{Corollary 8.} Let $\varphi$ be an Orlicz function.
\begin{enumerate}
\item If $\varphi\in \Delta_2$, then $Ces_\varphi=(Ces_\varphi)_a=C_\varphi$.
\item If $C : L^\varphi \rightarrow L^\varphi$ then $C_\varphi=Ces_\varphi$ if and only if $\varphi\in \Delta_2$.
\end{enumerate}
\bigskip
\textbf{Lemma 9.} If  $x\in C_\varphi$, then $\norm{x}_{Ces(\varphi)}=1$ if and only if $\rho_\varphi(x)=1$.
\begin{proof}
See the proof of Lemma 2.1 in [8].
\end{proof}
\bigskip
\noindent
\textbf{Theorem 10.} Let $\varphi$ be an Orlicz function satisfying $\varphi<\infty$.
\begin{enumerate}
\item If $I=[0,\infty)$, then the space $C_\varphi$ is strictly monotone if and only if $\varphi>0$.
\item If $I=[0,1]$ and $\lim_{u\rightarrow \infty} \varphi(u)/u= \infty$, then the space $C_\varphi$ is strictly monotone if and only if $\varphi>0$.
\end{enumerate}
\begin{proof}
\emph{Necessity}. Assume that $a_\varphi>0$. We have to consider two cases:

\noindent
(i). If $I=[0,\infty)$ then take $y_t=t\chi_{[0,1)}$  for $t\in \mathbb{R_+}$. We have
$$(Cy_t)(u)=t\chi_{[0,1)}(u)+\frac{t}{u} \chi_{[1,\infty)}(u).$$
We define the function
$$f(t)=\int_{0}^{\infty} \varphi\left(t\chi_{[0,1)}(u)+\frac{t}{u} \chi_{[1,\infty)}(u)\right)du=\rho_\varphi(y_t),$$
for $t>0$. Since $t/x \rightarrow 0$ as $x\rightarrow \infty$ and $a_\varphi>0$, so there exists $x_0=x_0(t)\in \mathbb{R_+}$ such that $t/x_0 \le a_\varphi$. Consequently, $f(t)\le \varphi(t) x_0$, which means that $f(t)$ takes finite values for $t>0$. Moreover, $f$ is convex function, so $f$ is continuous function on $\mathbb{R_+}$. This means that $f[\mathbb{R_+}]=\mathbb{R_+}$. In that case, from the Darboux property we find $\lambda>0$ satisfying $f(\lambda)=1$. Take $x_1>1$ with $\lambda/x_1 \le a_\varphi/2$. Let $z=a_\varphi \chi_{[x_1,x_1+1)}/2$. Then for any $x\ge x_1$ we have
$$C(y_\lambda+z)(x)=Cy_\lambda(x)+Cz(x)\le a_\varphi,$$
whence $\rho_\varphi(y_\lambda+z)=1$ and, from Lemma 9, we get $\norm{y_\lambda+z}_{Ces(\varphi)}=1$. Summing up, we built elements $u=y_\lambda$ i $v=y_\lambda+z$ such that $u\neq v$, $u\le v$ i $\norm{u}_{Ces(\varphi)}=\norm{v}_{Ces(\varphi}=1$. Moreover, $u,v\in C_\varphi$. This follows easily from the equality $C_\varphi=B_\varphi$ (see proof of Theorem 5). This means, that the space $C_\varphi$ isn't strictly monotone.

\noindent
(ii). Suppose $I=[0,1]$. For $a,b\in \mathbb{R}$ set $x=b\chi_{[0,a)}$. Then
$$(Cx)(t)=b\chi_{[0,a)}(t)+\frac{ab}{t}\chi_{[a,1)}(t).$$
Since $\lim_{u\rightarrow \infty} \varphi(u)/u = \infty$, there is a number  $b_0\in \mathbb{R}$ with $b_0>a_\varphi$ and $a_\varphi\varphi(b_0)/b_0 >1$. We define the function
$$f(a)=a\varphi(b_0)+\int_{a}^{\frac{ab_0}{a_\varphi}}\varphi\left(\frac{ab_0}{t}\right)dt,$$
for $a\in [0,a_\varphi/b_0 ]$. Clearly, $f(0)=0$ and $f(a_\varphi/b_0 )>1$. Moreover, we claim that $f$ is continuous on $[0,a_\varphi/b_0]$. Take $a_n\rightarrow a$, $a\in (0,a_\varphi/b_0]$ and put
$$A=\Big\{t\in \Big[0,\frac{a_\varphi}{b_0}\Big] : a\le t \le \frac{ab_0}{a_\varphi}\Big\},$$
$$A_n=\Big\{t\in \Big[0,\frac{a_\varphi}{b_0}\Big] : a_n\le t \le \frac{a_nb_0}{a_\varphi}\Big\},$$
for $n\in \mathbb{N}$. Then
$$0\le |f(a_n)-f(a)|\le |a_n-a|\varphi(b_0)+\Bigg|\int_{a}^{\frac{ab_0}{a_\varphi}} \varphi\left(\frac{ab_0}{t}\right)dt-\int_{a_n}^{\frac{a_nb_0}{a_\varphi}} \varphi\left(\frac{a_nb_0}{t}\right)dt\Bigg|$$
$$=|a_n-a|\varphi(b_0)+\Bigg| \int_{A\backslash A_n} \varphi\left(\frac{ab_0}{t}\right)dt
                             +\int_{A\cap A_n} \varphi\left(\frac{ab_0}{t}\right)-\varphi\left(\frac{a_nb_0}{t}\right)dt
                             -\int_{A_n\backslash A} \varphi\left(\frac{a_nb_0}{t}\right)dt\Bigg|$$
$$\le |a_n-a|\varphi(b_0)+\Bigg|\int_{A\backslash A_n} \varphi\left(\frac{ab_0}{t}\right)dt\Bigg|
                         +\Bigg|\int_{A\cap A_n} \varphi\left(\frac{ab_0}{t}\right)-\varphi\left(\frac{a_nb_0}{t}\right)dt\Bigg|
                         +\Bigg|\int_{A_n\backslash A} \varphi\left(\frac{a_nb_0}{t}\right)dt\Bigg|.$$
Now we have
$$0\le \Bigg|\int_{A\backslash A_n} \varphi\left(\frac{ab_0}{t}\right)dt\Bigg|\le m(A\backslash A_n)\sup\limits_{t\in A}\varphi(\frac{ab_0}{t})
=m(A\backslash A_n)\varphi(b_0)\rightarrow 0,$$
as $n\rightarrow \infty$. Similarly, we can show that $\Big|\int_{A_n\backslash A} \varphi\left(\frac{a_nb_0}{t}\right)dt \Big|\rightarrow 0$ as $n\rightarrow \infty$. Moreover, for $n$ large enough $a/2<a_n<2a$. Therefore, for $t\in A$
$$\Big|\varphi\left(\frac{ab_0}{t}\right)-\varphi\left(\frac{a_nb_0}{t}\right)\Big|\le 
\text{max}_{t_1,t_2\in [a_\varphi,2b_0]} |\varphi(t_1)-\varphi(t_2)|
=\eta(a,b_0)=\eta<\infty,$$
because $\varphi$ is continuous on the compact set $[a_\varphi,2b_0]$. Since $\big|\varphi(\frac{ab_0}{t})-\varphi(\frac{a_nb_0}{t})\big|\rightarrow 0$ pointwise on $A$, $\eta\chi_A$ is integrable majorant and $L^1 [0,1]$ is order continuous, so from Dominated Convergence Theorem we get
$$0\le \Bigg|\int_{A\cap A_n} \varphi\left(\frac{ab_0}{t}\right)-\varphi\left(\frac{a_nb_0}{t}\right)dt\Bigg|\le
\Bigg|\int_{A} \varphi\left(\frac{ab_0}{t}\right)-\varphi\left(\frac{a_nb_0}{t}\right)dt\Bigg|\rightarrow 0,$$
as $n\rightarrow \infty$. This means that $|f(a_n)-f(a)|\rightarrow 0$, so $f$ is continuous. This proves the claim.

By the Darboux property of $f$, there exists $a_1\in (0,a_\varphi/b_0 )$ such that $f(a_1 )=1$. Consider two elements
$$x_1=b_0\chi_{[0,a_1)},$$
and
$$x_2=b_0\chi_{[0,a_1)}+\left(a_\varphi-\frac{a_1b_0}{\delta}\right)\chi_{[\frac{\delta+1}{2},1)},$$
where $\delta>0$ is such that $a_1b_0/a_\varphi<\delta<1$. Obviously, $x_1\le x_2$ and $x_1\neq x_2$ and $x_1,x_2\in C_\varphi$. Furthermore, $\rho_\varphi(x_1)=\rho_\varphi (x_2)=1$. Indeed, we have the following inequalities
$$a_1<\frac{b_0}{a_\varphi}a_1<\delta<\frac{\delta+1}{2}<1,$$
so for $t\in [(\delta+1)/2,1)$
$$(Cx_2)(t)=\frac{a_1b_0}{t}+\left(a_\varphi-\frac{a_1b_0}{\delta}\right)\frac{1}{t}\int_{\frac{\delta+1}{2}}^{t}ds
<\frac{a_1b_0}{t}+a_\varphi-\frac{a_1b_0}{\delta}\le a_\varphi.$$
Thus, by the Lemma 9, $\norm{x_1}_{Ces(\varphi)}=\norm{x_2}_{Ces(\varphi)}=1$, whence the space $C_\varphi$ is not strictly monotone.

\noindent
\emph{Sufficiency.} Assume, that $a_\varphi=0$, $0\le x\le y$, $x\neq y$ and $x,y\in C_\varphi$. Without loss of generality, we can assume, that $\norm{x}_{Ces(\varphi)}=1$. From Lemma 9, $\rho_\varphi (x)=1$. We need only to show that $\rho_\varphi (y)>1$. Since $\varphi$ is superadditive, so for any  $x,y\in (C_\varphi)_{+}$  we have
$$\rho_\varphi(x+y)\ge \rho_\varphi(x)+\rho_\varphi(y).$$
Since $y-x\ge 0$, $y-x\neq 0$ and $\varphi >0$, we have $\rho_\varphi(y-x)>0$, and consequently
$$\rho_\varphi(y)=\rho_\varphi(x+(y-x))\ge \rho_\varphi(x)+\rho_\varphi(y-x)=1+\rho_\varphi(y-x)>1.$$
\end{proof}
\bigskip
\textbf{Theorem 11.} Let Let $\varphi$ be an Orlicz function with $\varphi<\infty$ and $C : L^\varphi\rightarrow L^\varphi$. Suppose additionally $\lim_{u\rightarrow \infty} \varphi(u)/u= \infty$ if $I=[0,1]$. The following conditions are equivalent:
\begin{enumerate}
\item $Ces_\varphi\in \text{(UM)}$,
\item $Ces_\varphi\in \text{(LLUM)}$,
\item $\varphi>0$ and $\varphi\in \Delta_2$.
\end{enumerate}
\begin{proof}
$\text{(i)}\Rightarrow \text{(ii)}$ by definition.

\noindent
$\text{(ii)}\Rightarrow \text{(iii)}$. If $a_\varphi>0$ then $C_\varphi \notin \text{(SM)}$ by Theorem 10. Moreover, see [12, Theorem 2.1], the following implication is true: if $X\in (H_l)$ then $X\in \text{(OC)}$. We can use the same proof to show that if $X\in\text{(LLUM)}$ then $X\in \text{(OC)}$. Therefore, if $\varphi\notin \Delta_2$ then $Ces_\varphi \notin \text{(OC)}$ by Theorem 7, whence $Ces_\varphi\notin \text{(LLUM)}$.

\noindent
$\text{(iii)}\Rightarrow \text{(i)}$. If $\varphi>0$ and $\varphi\in \Delta_2$ then $L_\varphi\in \text{(UM)}$, see e.g. [13, Theorem 7]. Consequently $Ces_\varphi\in \text{(UM)}$ by Fact.

\end{proof}

\section{ACKNOWLEDGEMENTS}

The second author (Pawe\l{} Kolwicz) is supported by the Ministry of Science and Higher Education of Poland, grant number 04/43/DSPB/0086.

\end{document}